\documentclass[preprint,1p,12pt]{elsarticle}

\makeatletter
\def\ps@pprintTitle{%
  \let\@oddhead\@empty
  \let\@evenhead\@empty
  \let\@evenfoot\@oddfoot
}
\makeatother

\usepackage{amssymb,amsmath,color}
\usepackage{graphicx}
\usepackage[ruled,vlined,linesnumbered]{algorithm2e}

\usepackage{url}
\urldef{\mailsa}\path|{alfred.hofmann, ursula.barth, ingrid.haas, frank.holzwarth,|
\urldef{\mailsb}\path|anna.kramer, leonie.kunz, christine.reiss, nicole.sator,|
\urldef{\mailsc}\path|erika.siebert-cole, peter.strasser, lncs}@springer.com|

\newcommand{\ot}{\leftarrow}

\newtheorem{theorem}{Theorem}
\newtheorem{corollary}[theorem]{Corollary}

\newdefinition{definition}[theorem]{Definition }
\newdefinition{remark}[theorem]{Remark }
\newdefinition{example}[theorem]{Example }
\newproof{proof}{Proof}

\begin{document}

\begin{frontmatter}

\title{On the solvability of bipolar max-product \\fuzzy relation equations with the product negation}

\author
{M. Eugenia Cornejo, David Lobo, Jes\'us Medina}

\address
{Department of Mathematics,
 University of  C\'adiz. Spain\\
\texttt{\{mariaeugenia.cornejo,david.lobo,jesus.medina\}@uca.es}}

\begin{abstract}

This paper studies the solvability of the max-product fuzzy relation equations in which a negation operator is considered. Specifically, the residuated negation of the product t-norm has been introduced in these equations in order to increase the flexibility of the standard fuzzy relation equations introduced by Sanchez in 1976. 
The solvability and the set of solutions of these bipolar equations have been studied in different scenarios, depending on the considered number of variables and equations. 
\end{abstract}

\end{frontmatter}

\section{Introduction}

Fuzzy relation equations (FREs) based on different max-t-norm compositions have become an extremely important research topic both from the theoretical and applicative perspective.  First of all, Sanchez introduced these equations by using the max-min composition for simulating the relationship between cause and effect in medical diagnosis problems~\cite{sanchez76,sanchez79}. Later, Pedrycz studied max-product fuzzy relation equations showing the applicative aspects of such equations in systems analysis, decision making and arithmetic of fuzzy numbers~\cite{PedryczGC:85}. Di Nola et al. provided a detailed study on max-continuous t-norm FREs and their applications to knowledge engineering in~\cite{Nola:1989}. Thenceforward, many researchers have continued investigating different procedures to solve FREs defined with either the max-min composition~\cite{chenwang07,chenwang02,Peeva:2004,Yeh2008}, the max-product composition~\cite{BOURKE1998,Loetamonphong:99,Markovskii,Peeva2007}, the max-Archimedean t-norm composition~\cite{LinJL,STAMOU2001,Wu2008} or other different compositions~\cite{Kohout80,Radim2012,baets99_eq,dm:mare,dm:ins2014,Medina2016,Ignjatovic20151,Medina2017:ija,Medina2017:ins402,Peeva2016,perfilieva08_fss}. 

As well as developing technical results on the resolution of FREs, the versatility of the applications of such equations in fuzzy sets theory has been increased. For instance, FREs play an important role in the treatment of images~\cite{Hirota1999, Hirota2002,Loia2005,Nobuhara2000,Nobuhara2002, Nobuhara2006}. Specifically, max-Lukasiewicz FREs are considered for images and videos compression and decompression  in~\cite{Loia2005}.  Fuzzy relation equations are also helpful to model the fuzzy inference systems considered in fuzzy control~\cite{Ciaramella2006146} and to represent restrictions in optimization problems~\cite{Li2008, Li2009,molai2010,yang2016}. In particular, max-product fuzzy relation equations are used for optimizing wireless communication management models in~\cite{yang2016}. FREs have also been related to other mathematical frameworks, such as fuzzy formal concept analysis~\cite{Cornejo2017,dm:mare,dm:ins2014}.

Many of the applications previously mentioned  require that the variables involved in fuzzy relation equations show a bipolar character. That is, these applications need equations containing unknown variables together with their logical negations simultaneously. 
For example, let us consider the equation $0.7*x=b$, where $*$ is the product t-norm, $x$ indicates the level of infection one patient has in the throat and $b$ the level of headache. This relation between infection and headache is modeled by the coefficient $0.7$.   If the doctor knows that the patient has an initial infection, with truth degree $0.5$, then we have that the patient has headache with truth degree: $0.35$. 
It is clear that this equation is a simple equation from the whole set of equations for simulating the considered knowledge system.
Suppose that the patient has not infection and has headache, with truth value $0.8$. In this case, the equation  $0.7*x=b$  does not work. However, if we consider  the product negation operator defined as $n(x)= 1$, if $x=0$, and $n(x)= 0$, otherwise, the previous equation can be conveniently reformulated obtaining the following one:
$$
(0.7*x)\vee (0.8*n(x))=b
$$ 
 which is solvable and captures the cases when the patient has the disease but does not have the symptoms. 
This new type of fuzzy relation equation, in which a negation operator is considered, is called bipolar fuzzy relation equation.  Notice that the use of the product negation increases the flexibility of the usual (unipolar) fuzzy relation equations and, at the same time, the complexity in their resolution, which becomes  even more difficult because the product negation is not involutive.

Bipolar fuzzy relation equations have been already employed usefully in optimization problems~\cite{Freson2013, Hu2016,Li2014,Liu2016, Zhou2016}. To the best of our knowledge, the specific literature on the resolution of bipolar fuzzy relation equations is really limited~\cite{Li2016, Liu2015}. Moreover, these papers are focused on the resolution of bipolar fuzzy relation equations based on max-min composition with the standard negation. Taking into account these considerations, we can ensure that the study on the solvability of bipolar FREs is an incipient research topic which has not been enough investigated so far.

This paper aims to provide under what conditions the solvability of bipolar max-product fuzzy relation equations with the product negation is guaranteed, continuing the initial study about this topic presented by the authors in~\cite{escim2017}. First of all, we will include  the  preliminary notions related to the calculus operators which will be used in bipolar max-product fuzzy relation equations.  In the following, we will introduce the notion of bipolar max-product FRE with the product negation and a unique unknown variable. Simple scenarios will be analyzed and a characterization on the solvability of such equations will be given. Providing  sufficient and necessary conditions to ensure when bipolar max-product FREs, with the product negation and different unknown variables, are solvable will be our next goal. Moreover, we will include different properties associated with the existence of the greatest/least solution or a finite number of maximal/minimal solutions for these last equations. We will also deal with the resolution of bipolar max-product FRE systems containing the product negation.  It is also convenient to mention that the results obtained in this paper differ markedly from those presented in~\cite{Li2016, Liu2015}, since the behaviour of the minimum t-norm and the product t-norm is clearly different. Some conclusions and prospects for future work are included at the end of the paper.

\section{Preliminary notions}

This section introduces basic notions and examples associated with the operators that we will use to make the computations in this paper.

Triangular norms are interesting operators which play an important role in different fields of mathematics such as probabilistic metric spaces~\cite{schweizer83,Serstnev}, decision making~\cite{Fodor1994,Grabisch1994}, statistics~\cite{Nelsen1999}, fuzzy sets theory and its applications~\cite{alsina83,Yager1980}. A detailed survey on the basic analytical and algebraic properties of triangular norms can be found in~\cite{Klement2004}. Below, we will include the formal definition of a triangular norm.

\begin{definition}\label{def:tnorm}
A binary operation $T\colon [0,1]\times [0,1]\to [0,1]$ is a \emph{triangular norm} (\emph{t-norm}) if the following properties are satisfied, for all $x,y,z\in [0,1]$:
\begin{itemize}
\item[(1)] $T(x,y)=T(y,x)$ \  (commutativity);
\item[(2)] If $x\leq y$ then $T(x,z)\leq T(y,z)$ \ (monotonicity);
\item[(3)] $T(x,1)=x$  \ (neutral element);
\item[(4)] $T(x,T(y,z))=T(T(x,y),z)$ \ (asociativity).
\end{itemize}
\end{definition}

T-norms are operators frequently used in applicative examples~\cite{pap,nguyen06,schweizer63} together with their residuated implications.

\begin{definition}\label{def:residimp}
Given a t-norm $T\colon [0,1]\times [0,1]\to [0,1]$, if there exists a binary operation $\ot_T \colon [0,1]\times [0,1]\to [0,1]$ order-preserving on the first argument and order-reversing on the second argument verifying the equivalence:
\begin{equation}\label{adjprop}
T(x,y)\leq z \ \hbox{ if and only if }  \  x\leq z\ot_T y
\end{equation}
then we say that $\ot_T$ is a \emph{residuated implication of \,$T$}. The pair $(T,\ot_T)$ is called \emph{residuated pair} and Equivalence~\eqref{adjprop} is called adjoint property.
\end{definition}

The residuated pairs of the most commonly used t-norms in the literature will be recalled in the following example.

\begin{example}\label{ex:tnorm}
Example of residuated pairs are the G\"odel, product and \L ukasiewicz  t-norms together with their corresponding residuated implications, which are defined on $[0,1]$, as you can see below:

$$
\begin{array}{lll}
T_{\text{P}} (x,y)= x \cdot y &  &   z\leftarrow_{\text{P}}x = \min\{1,z/x\} \\[2ex]
T_{\text{G}} (x,y)=\min\{x,y\} &  &   z\leftarrow_{\text{G}}x =  \begin{cases}1 & \hbox{if }  x\le z\\z&\hbox{otherwise}\end{cases} \\[3ex]
T_{\text{L}} (x,y)=\max\{0,x+y-1\} & & z\leftarrow_{\text{L}}x =\min\{1,1-x+z\}  \\[2ex]
\end{array}
$$
\qed
\end{example}

The residuated negation will be another fuzzy operator which will play a key role in this paper. These residuated negation operators are defined from the implication of a residuated pair. More information about this kind of negations can be found in~\cite{negationINS15}.
\begin{definition}\label{def:residneg}
Given an adjoint pair $(T,\ot_T)$ defined on $[0,1]$, the mapping $n_T\colon [0,1]\to [0,1]$ defined as $n_T(x)=0\leftarrow_T x$, for all $x\in[0,1]$, is called \emph{residuated negation}. We say that $n_T$ is a \emph{strong (or involutive) residuated negation} if the equality $x=n_T(n_T(x))$ holds for all $x\in[0,1]$.
\end{definition}

In the following example, we will introduce the residuated negations defined from the G\"odel, product and \L ukasiewicz  implications.

\begin{example}
Taking into account the definition of the operators introduced in Example~\ref{ex:tnorm}, it is easy to see that the residuated negation associated with the \L ukasiewicz implication is 
\begin{eqnarray*}
n_{\text{L}} (x)&=&1-x
\end{eqnarray*}
for all $x\in[0,1]$,  which is involutive and commonly known in the literature as the standard negation.  On the other hand, the G\"odel implication and the product implication coincide for all $x\in[0,1]$:
\begin{eqnarray*}
n_{\text{P}} (x)\, =\,\, n_{\text{G}} (x)&=& \begin{cases}1 & \hbox{if }  x=0 \\0&\hbox{otherwise}\end{cases}
\end{eqnarray*}
which are not involutive. This negation operator will be called \emph{product negation}.
\qed
\end{example}

In order to make the paper self-contained, we present the Intermediate Value Theorem~\cite{Bolzano1817} (english translation in~\cite{Russ1980}), which will be  used in different results later on.
\begin{theorem}[Intermediate Value Theorem]\label{th:int_value}
  Let $f\colon[a,b]\rightarrow \mathbb{R}$ be a continuous function and $c\in\mathbb{R}$ be any number between $f(a)$ and $f(b)$ inclusive. Then, there exists $x\in[a,b]$ such that $f(x)=c$.
\end{theorem}

From now on, we will focus on solving bipolar fuzzy relation equations based on the max-product t-norm composition. Our study will analyze the resolution of these equations when the considered negation operator is the non-involutive product negation.

\section{Bipolar max-product FREs with the product negation}\label{sec:bmpn}

This section introduces an interesting study on how to solve bipolar max-product fuzzy relation equations with the product negation. The results will be progressively presented according to their complexity, that is: (1) we will firstly analyze when a bipolar max-product FRE with the product negation and a unique unknown variable has solution; (2) hereunder, we will show under what conditions bipolar max-product FREs with different unknown variables are solvable and have a greatest (least, respectively) solution or a finite number of maximal (minimal, respectively) solutions; (3) we will finish this section with a characterization on the solvability of bipolar max-product FRE systems.

\subsection{Solving bipolar max-product FREs with one unknown variable}
The most simple scenario which can be considered in this paper is the study of a bipolar max-product FRE containing a unique unknown variable. In order to provide a characterization on the solvability of these equations, we need to include the following formal definition.

\begin{definition}\label{def:bmpp}
Let $a^+,a^-,b\in[0,1]$, $x$ be an unknown variable in $[0,1]$, $\ast$ be the product t-norm, $\vee$ be the maximum operator and $n_P$ be the product negation. Equation~\eqref{eq:adj_simple_0} is called \emph{bipolar max-product fuzzy relation equation with the product negation}.
\begin{equation}\label{eq:adj_simple_0}
    (a^+ \ast x) \vee (a^-\ast n_P (x))=b
  \end{equation}
\end{definition}

Before presenting under what conditions Equation~\eqref{eq:adj_simple_0} is solvable, we will analyze its solvability when either $a^+$, $a^-$ or $b$ is equal to $0$. {The following result characterizes the solvability of Equation~\eqref{eq:adj_simple_0} in the previous three different cases.

\begin{theorem}\label{th:adj_simple}
Let $a^+,a^-,b\in[0,1]$ and $x$ be an unknown variable belonging to $[0,1]$. Then, the following statements hold:
\begin{itemize}
\item [(1)] If $a^+=0$, then Equation~\eqref{eq:adj_simple_0} is solvable if and only if $b=0$ or $a^-=b$.
\item [(2)] If $a^-=0$, then Equation~\eqref{eq:adj_simple_0} is solvable if and only if $b\leq a^+$.
\item [(3)] If $b=0$, then Equation~\eqref{eq:adj_simple_0} is solvable if and only if $a^+=0$ or $a^-=0$.
\end{itemize}
\begin{proof}
    In order to demonstrate Statement (1), suppose that $a^+=0$. In that case, we obtain that  Equation~\eqref{eq:adj_simple_0} becomes into
    $a^-\ast n_P (x)=b$. Taking into account the definition of the product negation, we can obtain that $b=0$ if and only if any $x\in(0,1]$ is a solution of Equation~\eqref{eq:adj_simple_0}. In addition, the equality $a^-=b$ holds if and only if $x=0$ is a solution Equation~\eqref{eq:adj_simple_0}. Therefore, Statement (1) is verified.

    Now, assume that $a^-=0$. As a consequence, Equation~\eqref{eq:adj_simple_0} is given by $a^+\ast x=b$. Due to the product t-norm is a continuous order-preserving mapping on $[0,1]$ and the equalities $a^+\ast 0=0$ and $a^+\ast 1=a^+$ hold, we can assert by the Bolzano's theorem that there exists $x\in[0,1]$ such that $a^+\ast x=b$ if and only if $b\leq a^+$. Hence, we obtain that Statement (2) is satisfied.

    Finally, if $b=0$, then Equation~\eqref{eq:adj_simple_0} becomes into $(a^+ \ast x) \vee (a^-\ast n_P (x))=0$, which is solvable if and only if $a^+ \ast x=0$ and $a^-\ast n_P (x)=0$. According to the definition of the product negation, we can deduce that Equation~\eqref{eq:adj_simple_0} is solvable if and only if $a^+=0$ or $a^-=0$. As a result, we conclude that Statement (3) holds.
  \qed
\end{proof}
\end{theorem}
}

After studying the most trivial cases which can be given in a bipolar max-product FRE, we will continue our research assuming that
each known variable appearing in bipolar max-product FREs is different from zero.

A characterization on the solvability of bipolar max-product FREs with the product negation and a unique unknown variable is provided by the next theorem.

\begin{theorem}\label{th:adj_simple}
Let $a^+,a^-,b\in(0,1]$ and $x$ be an unknown variable belonging to $[0,1]$. The bipolar max-product FRE given by Equation~\eqref{eq:adj_simple_0} is solvable if and only if \,$a^-=b$ or $b\leq a^+$. In this case, at most two solutions exists: $x=0$ or/and $b/a^+$, which  are related to the conditions $a^-=b$ and $b\leq a^+$, respectively.
\end{theorem}

\begin{proof}
First of all, we will prove that if $a^-=b$ or $a^+\geq b$ then Equation~\eqref{eq:adj_simple_0} has solution.  On the one hand, if $a^-=b$ then $x=0$ is clearly a solution of Equation~\eqref{eq:adj_simple_0}. On the other hand, we suppose that \,$a^+\geq b$. Now, we consider the mapping $f\colon [0,1]\to[0,1]$ defined as $f(x)=a^+*x$, which is continuous. Since $b\in (0,1]$ is a number between $f(0)=a^+ * 0=0$ and $f(1)=a^+ * 1=a^+$, applying the Intermediate Value Theorem (Theorem~\ref{th:int_value}), there exists $y \in[0,1]$ such that $f(y)=a^+ * y=b$. In addition, we can ensure that the previous value $y \in[0,1]$ is unique since $f$ is a strictly order-preserving mapping. Specifically, by the definition of the product implication,  $y=b\leftarrow_P a^+$. In the following, we will demonstrate that $y=b\leftarrow_P a^+$ is a solution of Equation~\eqref{eq:adj_simple_0}:
\begin{eqnarray*}
(a^+ \ast y) \vee (a^-\ast n_P (y))&=& (a^+ \ast (b\leftarrow_P a^+)) \vee (a^-\ast n_P (b\leftarrow_P a^+))\\
&=&(a^+ \ast b/a^+) \vee (a^-\ast n_P (b/a^+))\\
&=& b \vee (a^-\ast 0)=b\vee0=b
\end{eqnarray*}
being $b\neq0$ by hypothesis. Consequently, we obtain that Equation~\eqref{eq:adj_simple_0} is solvable.

In order to prove the counterpart, we will suppose that $a^-\neq b$ and $a^+< b$, and we will prove that Equation~\eqref{eq:adj_simple_0} is not solvable. Clearly, for each $x\in[0,1]$, we have that the inequality $a^+\ast x<b$ holds. Therefore, Equation~\eqref{eq:adj_simple_0} is solvable if and only if there exists $x\in[0,1]$ such that $a^-\ast n_P(x)=b$. However, this last condition is not verified, since $n_P$ can only take either the value $0$ or the value $1$ and, by hypothesis, $a^-\neq b$ and $b\neq 0$. Hence, we can conclude that Equation~\eqref{eq:adj_simple_0} is not solvable.

Finally, when Equation~\eqref{eq:adj_simple_0} is solvable, it has at most two solutions which are given by $x=0$ and $x=b\leftarrow_P a^+$. It is easy to see that $x=b\leftarrow_P a^+=b/a^+$ is the greatest solution when both $x=0$ and $x=b\leftarrow_P a^+$ solve Equation~\eqref{eq:adj_simple_0}.
\qed
\end{proof}

An illustrative example will be shown in order to clarify the previous theorem.

\begin{example}\label{ex:adj_simple}
A simple bipolar max-product  FRE is provided by Equation~\eqref{eq:ex:adj_simple} and we will apply Theorem~\ref{th:adj_simple} to check whether such equation is solvable or not.
\begin{equation}\label{eq:ex:adj_simple}
(0.5 \ast x) \vee (0.2\ast n_P (x))=0.3
\end{equation}
Observe that, $0.2\neq0.3$ and $0.5\geq 0.3$. Therefore, we can ensure that the hypothesis required in Theorem~\ref{th:adj_simple} are satisfied and Equation~\eqref{eq:ex:adj_simple} is solvable. Indeed, $x=0.3\leftarrow_P 0.5=0.6$ is a solution, as the next equality shows:
$$
(0.5 \ast 0.6) \vee (0.2\ast n_P (0.6))=(0.5 \ast 0.6) \vee (0.2\ast 0)=0.3 \vee 0=0.3
$$
Moreover, $x=0.3\leftarrow_P 0.5=0.6$ is the unique solution since $x=0$ does not verify Equation~\eqref{eq:ex:adj_simple}, as we can see below:
$$
(0.5 \ast 0) \vee (0.2\ast n_P (0))=(0.5 \ast 0) \vee (0.2\ast 1)=0 \vee 0.2=0.2\neq0.3
$$
\qed
\end{example}

Once we have studied bipolar max-product fuzzy relation equations with the product negation and containing only one unknown variable, our following goal consists in increasing the number of unknown variables appearing in a bipolar max-product FRE with the product negation and studying sufficient and necessary conditions to guarantee its solvability.

\subsection{Solving bipolar max-product FREs with different unknown variables}
This section includes results associated with the solvability of bipolar max-product FREs with the product negation and a finite number of different unknown variables. Besides, we will present the conditions under which the existence of either a greatest (least, respectively) solution or a finite number of maximal (minimal, respectively) solutions can be guaranteed in these equations.

{To begin with this section, we will present a result which characterizes the solvability of a bipolar max-product FRE when the independent term is equal to zero. The consideration of this case will allow us to obtain under what conditions a general bipolar max-product FRE with different unknown variables is solvable in an easier way.

\begin{theorem}\label{th:adj_bzero}
Given $a_j^+,a_j^-\in[0,1]$ and $x_j$ an unknown variable belonging to $[0,1]$, for all $j\in\{1,\dots,m\}$. The bipolar max-product  fuzzy relation equation
  \begin{equation}\label{eq:adjprod2_bzero}
    \bigvee_{j=1}^m(a_j^+ \ast x_j) \vee (a_j^-\ast n_P (x_j))=0
  \end{equation}
  is solvable if and only if, for each $j\in\{1,\dots,m\}$, either the equality $a_j^+=0$ or $a_j^-=0$ is satisfied.
\end{theorem}
\begin{proof}
  First of all, we will suppose that, for each $j\in\{1,\dots,m\}$, either the equality $a_j^+=0$ or $a_j^-=0$ is satisfied, and we will find a solution of Equation~\eqref{eq:adjprod2_bzero}. Given $k\in\{1,\dots,m\}$, if $a_k^+=0$, then the variable  $\hat{x}_k$ of the solution will be $1$, since, in this case, we obtain:
  \[(a_k^+ \ast \hat{x}_k) \vee (a_k^-\ast n_P (\hat{x}_k))=0\vee (a_k^-\ast 0)=0\]
  Otherwise, if $a_k^-=0$, then we consider $\hat{x}_k=0$, which leads us to the same result:
  \[
  (a_k^+ \ast \hat{x}_k) \vee (a_k^-\ast n_P (\hat{x}_k))=(a_k^-\ast 0)\vee0=0
  \]
Therefore, the obtained tuple $(\hat{x}_1,\dots,\hat{x}_m)$ verifies that: 
  \[
  \bigvee_{j=1}^m(a_j^+ \ast \hat{x}_j) \vee (a_j^-\ast n_P (\hat{x}_j))=0
  \]
  and, as a consequence, it  is a solution of Equation~\eqref{eq:adjprod2_bzero}.

  Now, we suppose that Equation~\eqref{eq:adjprod2_bzero} is solvable and the tuple $(x_1,\dots,x_m)$ is a solution of such equation. We will assume that there exists $k\in\{1,\dots,m\}$ such that $a_k^+>0$ and $a_k^->0$, and we will obtain a contradiction.

If $x_k=0$, then $n_P(x_k)=1$, and therefore the chain of inequalities below holds
  \[\bigvee_{j=1}^m(a_j^+ \ast x_j) \vee (a_j^-\ast n_P (x_j))\geq (a_k^-\ast n_P (x_k))=a_k^->0\]
  which contradicts that   $(x_1,\dots,x_m)$ is  a solution of Equation~\eqref{eq:adjprod2_bzero}.
    On the other hand, if $x_k>0$, due to the fact that $a_k^+>0$, we obtain the same strict inequality:
  \[
  \bigvee_{j=1}^m(a_j^+ \ast x_j) \vee (a_j^-\ast n_P (x_j))\geq (a_k^+\ast x_k)>0
  \]
  which also contradicts the hypothesis. Thus,  we can ensure that either the equality $a_j^+=0$ or $a_j^-=0$ holds, for all $j\in\{1,\dots,m\}$.
  \qed
\end{proof}
}
The following theorem characterizes the  solvability of the proper bipolar max-product FREs containing different unknown variables.

\begin{theorem}\label{th:adj}
Given $a_j^+,a_j^-\in[0,1]$, $b\in(0,1]$ and $x_j$ an unknown variable belonging to $[0,1]$, for all $j\in\{1,\dots,m\}$. The bipolar max-product  fuzzy relation equation
  \begin{equation}\label{eq:adjprod2}
    \bigvee_{j=1}^m(a_j^+ \ast x_j) \vee (a_j^-\ast n_P (x_j))=b
  \end{equation}
  is solvable if and only if $b\leq\max\{a_j^+\!\!\mid\! j\!\in\!\{1,\dots,m\}\}$ or there exists $k\in\{1,\dots,m\}$ such that $a_k^-=b$.
\end{theorem}
\begin{proof}
We will suppose that $b\leq\max\{a_j^+\mid j\in\{1,\dots,m\}\}$ and we will prove that Equation~\eqref{eq:adjprod2} is solvable. Clearly, $J=\{j\in\{1,\dots,m\}\mid b\leq a_j^+\}$ is a non-empty set. Furthermore, by definition of the product residuated implication, the following statements hold, for each $j\in\{1,\dots,m\}$:
\begin{itemize}
\item[(i)] If $j\in J$, then $b\leftarrow_P a_j^+=b/a_j^+$ \footnote{Notice that, if $a_j^+=0$, then $a_j^+<b$, and thus $a_j^+\notin J$.} and we obtain the following equalities $a_j^+*(b\leftarrow_P a_j^+)=b$ and $a_j^-\ast n_P (b\leftarrow_P a_j^+)=0$.
\item[(ii)] If $j\notin J$, then we obtain that  $b\leftarrow_P a_j^+=1$ and the following expressions hold $a_j^+*(b\leftarrow_P a_j^+)=a_j^+<b$ and $a_j^-\ast n_P (b\leftarrow_P a_j^+)=0$.
\end{itemize}
From (i) and (ii), we can conclude that the tuple $(b\leftarrow_P a_1^+,\dots,b\leftarrow_P a_m^+)$ is a solution of Equation~\eqref{eq:adjprod2}.

Now, we will suppose that there exists $k\in\{1,\dots,m\}$ such that $a_k^-=b$ and $\max\{a_j^+\mid j\in\{1,\dots,m\}\}< b$ and we will demonstrate that Equation~\eqref{eq:adjprod2} is solvable. Notice that, if $b\leq\max\{a_j^+\mid j\in\{1,\dots,m\}\}$, then we are in the case above, and thus Equation~\eqref{eq:adjprod2} is solvable.
By hypothesis, $a_j^+<b$ for each $j\in\{1,\dots,m\}$, and therefore, $b\leftarrow_P a_j^+=1$. As a consequence, $a_j^+*(b\leftarrow_P a_j^+)=a_j^+<b$ and $a_j^-\ast n_P (b\leftarrow_P a_j^+)=0$, for each $j\in\{1,\dots,m\}$. Hence, we obtain straightforwardly that the tuple $(\hat{x}_1,\dots,\hat{x}_k,\dots,\hat{x}_m)$ with $\hat{x}_k=0$ and $\hat{x}_j=1$, for each $j\in\{1,\dots,m\}$ such that $j\neq k$, is  a solution of Equation~\eqref{eq:adjprod2}.

The proof of the counterpart is analogous to the proof of Theorem~\ref{th:adj_simple}.
  \qed
\end{proof}

From now on, explanatory examples will be included for a better understanding of each result.

\begin{example}\label{ex:adj}
This example will consider the bipolar max-product FRE with three unknown variables $x_1,x_2,x_3\in[0,1]$ given by the following equation:

{ \footnotesize \begin{equation}\label{eq:ex:adj}
(0.4 \ast x_1) \vee (0.7\ast n_P (x_1))\vee(0.2 \ast x_2) \vee (0.1\ast n_P (x_2))\vee(0.5 \ast x_3) \vee (0.2\ast n_P (x_3))=0.3
 \end{equation} }

We can ensure that Equation~\eqref{eq:ex:adj} is solvable because the hypothesis required in Theorem~\ref{th:adj} are verified, that is, $\max\{0.4,0.2,0.5\}=0.5\geq0.3$. Following an analogous reasoning to the one given in the proof of such theorem, we obtain that the tuple $(0.75,1,0.6)$ is a solution of Equation~\eqref{eq:ex:adj}. In addition, we can make simple computations in order to obtain all possible solutions of Equation~\eqref{eq:ex:adj}. For instance, the tuple $(0.5,0,0.6)$ is another solution of Equation~\eqref{eq:ex:adj}:

{ \footnotesize
\begin{eqnarray*}
(0.4 \ast 0.5) \vee (0.7\ast n_P (0.5))\vee(0.2 \ast 0) \vee (0.1\ast n_P (0))\vee(0.5 \ast 0.6) \vee (0.2\ast n_P (0.6)) &=& \\
 0.2 \vee 0\vee0 \vee 0.1\vee0.3 \vee 0 &=&\\ 0.3
\end{eqnarray*}
}\qed
\end{example}

In view of the results obtained in the previous example, one can ask whether $(0.75,1,0.6)$ is the greatest solution of Equation~\eqref{eq:ex:adj} or if there exists a finite number of maximal solutions for such equation. A similar question with respect to the existence of its least solution or different minimal solutions can also arise. These issues will be analyzed in the following section.

\subsubsection{Computing maximal and minimal solutions}

Next, we are interested in knowing when a bipolar max-product FRE with the product negation and different variables has either a greatest (least, respectively) solution or a finite number of maximal (minimal, respectively) solutions.  To reach this goal, we need to introduce the following results.

\begin{theorem}\label{th:adj_max}
Given $a_j^+,a_j^-\in[0,1]$, $b\in(0,1]$ and $x_j$ an unknown variable in $[0,1]$, for each $j\in\{1,\dots,m\}$.  If Equation~\eqref{eq:adjprod2} is a solvable bipolar max-product  FRE, then the following statements hold:
\begin{enumerate}
\item[(1)] If $b\leq\max\{a_j^+\mid j\in\{1,\dots,m\}\}$, the set of solutions of Equation~\eqref{eq:adjprod2} has a greatest element.
\item[(2)] If $a_j^+< b$ for each $j\in\{1,\dots,m\}$, then the set of maximal solutions of Equation~\eqref{eq:adjprod2} is finite. Moreover, the number of maximal solutions is:
        \[card(\{k\in\{1,\dots,m\}\mid a_k^-=b\})\]
  \end{enumerate}
\end{theorem}

\begin{proof}\
In order to prove Statement (1), suppose that the inequality $b\leq\max\{a_j^+\mid j\in\{1,\dots,m\}\}$ holds. By an analogous reasoning to the proof given in Theorem~\ref{th:adj}, we obtain that the tuple $(b\leftarrow_P a_1^+,\dots,b\leftarrow_P a_m^+)$ is a solution of Equation~\eqref{eq:adjprod2}. In the following, we will prove that $(b\leftarrow_P a_1^+,\dots,b\leftarrow_P a_m^+)$ is the greatest solution of Equation~\eqref{eq:adjprod2}, by reduction to the absurd.

Hence, we suppose that there exists a tuple $(x_1,\dots,x_m)$ being a solution of Equation~\eqref{eq:adjprod2} such that $(b\leftarrow_P a_1^+,\dots,b\leftarrow_P a_m^+)<(x_1,\dots,x_m)$. Then, there exists $k\in\{1,\dots,m\}$ such that $x_k>b\leftarrow_P a_k^+$. According to the adjoint property, we can assert that $a_k^+*x_k> b$. As a result, we obtain that

\begin{equation*}
\bigvee_{j=1}^m(a_j^+ \ast x_j) \vee (a_j^-\ast n_P (x_j))=\bigvee_{j=1}^m(a_j^+ \ast x_j)\geq a_k^+*x_k>b
\end{equation*}

which lead us to conclude that $(x_1,\dots,x_m)$ is not a solution of Equation~\eqref{eq:adjprod2}, in contradiction with the hypothesis. Hence, we can ensure that $(b\leftarrow_P a_1^+,\dots,b\leftarrow_P a_m^+)$ is the greatest solution of Equation~\eqref{eq:adjprod2}. That is, Statement (1) holds.

Now, assume that $a_j^+< b$ for each $j\in\{1,\dots,m\}$, in order to demonstrate Statement (2). Given $A=\{k\in\{1,\dots,m\}\mid a_k^-=b\}$ and $B$ the set of maximal solutions of Equation~\eqref{eq:adjprod2}, we consider the mapping $f\colon A\to B$ defined as $f(k)=(1,\dots,1,0,1,\dots,1)$, being the element $0$ in the $k$-th position of the tuple. In order to prove that the cardinal of $A$ coincides with the cardinal of $B$, we will demonstrate that $f$ is a bijection between $A$ and~$B$.

First of all, let us see that $f$ is well-defined. Consider a fixed $k\in A$. By an analogous reasoning to the proof of Theorem~\ref{th:adj}, we can ensure that the tuple $(1,\dots,1,0,1,\dots,1)$ is a solution of Equation~\eqref{eq:adjprod2}, being the element $0$ in the $k$-th position of the tuple. Suppose that there exists another solution $({\hat{x}_1},\dots,{\hat{x}_k},\dots,{\hat{x}_m})$ of  Equation~\eqref{eq:adjprod2}  such that  $(1,\dots,1,0,1,\dots,1)<({\hat{x}_1},\dots,{\hat{x}_k},\dots,{\hat{x}_m})$. Therefore, ${\hat{x}_k}>0$ and ${\hat{x}_j}=1$, for each $j\in\{1,\dots,m\}$ with $j\neq k$, in order to $({\hat{x}_1},\dots,{\hat{x}_k},\dots,{\hat{x}_m})$ be a greater solution. As a consequence, $n_P({\hat{x}_j})= 0$ for each $j\in\{1,\dots,m\}$ and so,
\begin{eqnarray*}
\bigvee_{j=1}^m(a_j^+ \ast x_j) \vee (a_j^-\ast n_P (x_j))&=&\bigvee_{j=1}^m(a_j^+ \ast x_j)\\
&=&a_1^+\vee\dots\vee (a_k^+*x_k)\vee\dots\vee a_m^+\\
&<& b
\end{eqnarray*}
since, by hypothesis, $a_j^+< b$ for each $j\in\{1,\dots,m\}$. Consequently, the tuple $({\hat{x}_1},\dots,{\hat{x}_k},\dots,{\hat{x}_m})$ is not a solution of  Equation~\eqref{eq:adjprod2}. Thus, $(1,\dots,1,0,1,\dots,1)$ is a maximal solution of Equation~\eqref{eq:adjprod2}, which permits us to assert that $f$ is well-defined. Moreover, given $k_1,k_2\in A$ with $k_1\neq k_2$, then $f(k_1)$ is clearly different than $f(k_2)$, since the element $0$ is in the \emph{$k_1$}-th position in $f(k_1)$ while the element $0$ appears in the \emph{$k_2$}-th position in $f(k_2)$. Therefore, $f$ is an order-embedding mapping.

Now, we will prove that $f$ is also onto. Let $(x_1,\dots,x_m)\in B$. Clearly, as $a_j^+< b$ for each $j\in\{1,\dots,m\}$, we have that $a_j^+*x_j< b$. Therefore, as the tuple $(x_1,\dots,x_m)$ is a solution of Equation~\eqref{eq:adjprod2}, we can ensure that there exists $k\in\{1,\dots,m\}$ such that $a_k^-*n_P(x_k)=b$. According to the definition of the product negation, we obtain that $n_P(x_k)$ can only take the values $0$ and $1$. Since $b\neq 0$, we deduce that $n_P(x_k)=1$, and as a consequence, $x_k=0$. This fact allows us to ensure that $a_k^-=b$, that is, $k\in A$. Due to $(x_1,\dots,x_m)$ is a maximal solution of Equation~\eqref{eq:adjprod2} and the inequality $a_j^+*x_j< b$ holds, for each $j\in\{1,\dots,m\}$, we can assert that $x_j=1$ for each $j\in\{1,\dots,m\}$ with $j\neq k$. Hence, we obtain that $f(k)=(1,\dots,1,0,1,\dots,1)=(x_1,\dots,x_m)$, being the element $0$ in the $k$-th position of the tuple.

Due to the mapping $f$ is a bijection, we can conclude that the number of maximal solutions is $card(\{k\in\{1,\dots,m\}\mid a_k^-=b\})$. Equivalently, Statement (2) is satisfied.
  \qed
\end{proof}

As a  consequence of Theorem~\ref{th:adj_max}, the next corollary arises.

\begin{corollary}
Given $a_j^+,a_j^-\in[0,1]$, $b\in(0,1]$ and $x_j$ an unknown variable in $[0,1]$, for each $j\in\{1,\dots,m\}$.  Consider that Equation~\eqref{eq:adjprod2} is a solvable bipolar max-product  FRE, then the following statements hold:
\begin{enumerate}
\item[(1)] If $\max\{a_j^+\mid j\in\{1,\dots,m\}\}\geq b$, then the greatest solution of Equation~\eqref{eq:adjprod2} is given by the tuple $(b\leftarrow_P a_1^+,\dots,b\leftarrow_P a_m^+)$.
\item[(2)] If $a_j^+< b$ for each $j\in\{1,\dots,m\}$, then the set of maximal solutions of Equation~\eqref{eq:adjprod2} is given by:
    \[\{(1,\dots,1, x_k,1,\dots,1)\mid  x_k=0 \hbox{ with } k\in K^-_P\}\]
where $K^-_P= \{k\in\{1,\dots,m\}\mid a_k^-=b$\}
\end{enumerate}
\end{corollary}

We also have to distinguish different cases to ensure the existence of  the least solution and the set of minimal solutions of a solvable bipolar max-product  FRE with different variables. It is important to mention that the set of minimal solutions of such equation can be empty, as the following theorem shows.

\begin{theorem}\label{th:adj_min}
Given $a_j^+,a_j^-\in[0,1]$, $b\in(0,1]$ and $x_j$ an unknown variable belonging to  $[0,1]$, for each $j\in\{1,\dots,m\}$.   Consider that Equation~\eqref{eq:adjprod2} is a solvable bipolar max-product  FRE.
  Then, the following statements hold:
  \begin{enumerate}
    \item[(1)] If there exists $k\in\{1,\dots,m\}$ such that $a_k^-=b$ and $a_j^-\leq b$, for each $j\in\{1,\dots,m\}$, then the set of solutions of Equation~\eqref{eq:adjprod2} has a least element.
    \item[(2)] If there exist $k_1,k_2\in\{1,\dots,m\}$ such that $a_{k_1}^-=b$ and $a_{k_2}^->b$, then the set of solutions of Equation~\eqref{eq:adjprod2} has no minimal elements.
    \item [(3)] If $a_j^-\neq b$ for each $j\in\{1,\dots,m\}$, then the set of minimal solutions of Equation~\eqref{eq:adjprod2} is finite. Moreover, the number of minimal solutions is:
        \[card(\{k\in\{1,\dots,m\}\mid a_k^+\geq b\text{ and } a_j^- < b \text{ for each }j\neq k\})\]
  \end{enumerate}
\end{theorem}
\begin{proof}
First of all, to prove Statement (1), we will suppose that there exists $k\in\{1,\dots,m\}$ such that $a_k^-=b$ and $a_j^-\leq b$, for each $j\in\{1,\dots,m\}$. Then, the tuple  $(0,\dots,0)$ is clearly a solution of Equation~\eqref{eq:adjprod2}, since:
\[\bigvee_{j=1}^m(a_j^+ \ast 0) \vee (a_j^-\ast n_P(0))=\bigvee_{j=1}^m(a_j^-\ast n_P(0))=\bigvee_{j=1}^m a_j^-=b\]
Clearly, any other solution of Equation~\eqref{eq:adjprod2} is greater than $(0,\dots,0)$. As a consequence, $(0,\dots,0)$ is the least element of the set of solutions of Equation~\eqref{eq:adjprod2}. That is, Statement (1) is satisfied.

Now, in order to demonstrate Statement (2), we suppose that there exist $k_1,k_2\in\{1,\dots,m\}$ such that $a_{k_1}^-=b$ and $a_{k_2}^->b$. We will prove that the set of solutions of Equation~\eqref{eq:adjprod2} has no minimal elements by reduction to the absurd.

Let us assume that $(x_1,\dots,x_m)$ is a minimal solution of Equation~\eqref{eq:adjprod2}. First of all, we will see that $x_{k_1}=0$ and $x_{k_2}>0$. On the one hand, suppose that $x_{k_1}>0$. Since equalities $a_{k_1}^+*0=0$ and $a_{k_1}^-*n_P(0)=a_{k_1}^-=b$ are straightforwardly satisfied, we obtain that the tuple $(x_1,\dots,x_{k_1-1},0,x_{k_1+1},\dots, x_m)$ is another solution of Equation~\eqref{eq:adjprod2}, which is strictly smaller than $(x_1,\dots,x_m)$. Thus, we can ensure than $x_{k_1}=0$. On the other hand, let us assume that $x_{k_2}=0$. In this case, as $a_{k_2}^- >b$, the following inequality holds
\[\bigvee_{j=1}^m(a_j^+ \ast x_j) \vee (a_j^-\ast n_P (x_j))\geq a_{k_2}^-* n_P(x_{k_2})=a_{k_2}^->b\]
As a result, the tuple $(x_1,\dots,x_m)$ is not a solution of Equation~\eqref{eq:adjprod2}. Consequently, we can ensure that $x_{k_2}>0$.
Therefore, we can assert that $x_{k_1}=0$ and $x_{k_2}>0$.
In addition, notice that $a_{k_1}^-*n_P(x_{k_1})=b$. Due to the fact that $(x_1,\dots,x_m)$ is a solution of Equation~\eqref{eq:adjprod2}, the inequality $a_{k_2}^+*x_{k_2}\leq b$ is verified. Then, as the product t-norm is a monotonic operator, we obtain that inequality $a_{k_2}^+*\frac{x_{k_2}}{2}\leq b$ also holds. Notice that, since $x_{k_2}>0$, then $\frac{x_{k_2}}{2}<x_{k_2}$. Since $a_{k_1}^-*n_P(x_{k_1})=b$, the tuple $(x_1,\dots,x_{k_2-1},\frac{x_{k_2}}{2},x_{k_2+1}\dots,x_m)$ is a solution of Equation~\eqref{eq:adjprod2}, which satisfies $(x_1,\dots,x_{k_2-1},\frac{x_{k_2}}{2},x_{k_2+1}\dots,x_m)<(x_1,\dots,x_m)$.
As a result, we conclude that $(x_1,\dots,x_m)$ is not a minimal solution, in contradiction with the hypothesis, and thus we can assert that Statement (2) holds.

Finally, we will show that Statement (3) is verified. Given $A=\{k\in\{1,\dots,m\}\mid a_k^+\geq b\text{ and }a_j^- \leq b \text{ for each }j\neq k\}$ and $B$ the set of minimal solutions of Equation~\eqref{eq:adjprod2}, we will demonstrate that the mapping $f\colon A\to B$, which associates each $k\in A$ with a tuple $(0,\dots,0,b\leftarrow_P a_k^+,0\dots,0)$, being the element $b\leftarrow_P a_k^+$ in the \emph{k}-th position of the tuple, is a bijection.

Now, we will prove that $f$ is well-defined. Consider a fixed $k\in A$. Since $a_k^+\geq b$ and $b\neq 0$, we obtain that $b\leftarrow_P a_k^+= b/a_k^+>0$ and consequently, the following equalities $a_k^+*(b\leftarrow_P a_k^+)=b$ and $a_k^-*n_P(b\leftarrow_P a_k^+)=0$ are satisfied. In addition, since the inequality $a_j^- <b$ holds for each $j\in\{1,\dots,m\}$ with $j\neq k$, we can ensure that the tuple $(0,\dots,0,b\leftarrow_P a_k^+,0,\dots,0)$ is a solution of Equation~\eqref{eq:adjprod2}. Suppose that there exists another different tuple $(x_1,\dots,x_k,\dots,x_m)$ being solution of Equation~\eqref{eq:adjprod2} and satisfying that $(x_1,\dots,x_k,\dots,x_m)<(0,\dots,0,b\leftarrow_P a_k^+,0,\dots,0)$. Then, we have that $x_k<b\leftarrow_P a_k^+$ and $x_j=0$, for each $j\in\{1,\dots,m\}$ with $j\neq k$. Clearly, as the product t-norm is strictly order-preserving, the inequality $a_k^+ * x_k<b$ is also verified. Finally since, by hypothesis, $a_j^-\neq b$ for each $j\in\{1,\dots,m\}$, we can conclude that $a_j^-*n_P(x_j)\neq b$ for each $j\in\{1,\dots,m\}$, and thus, $(x_1,\dots,x_m)$ is not a solution of Equation~\eqref{eq:adjprod2}. Therefore, $(0,\dots,0,b\leftarrow_P a_k^+,0,\dots,0)$ is a minimal solution of Equation~\eqref{eq:adjprod2} and so, we can ensure that $f(k)=(0,\dots,0,b\leftarrow_P a_k^+,0,\dots,0)\in B$, that is, the mapping $f$ is well-defined.

In the following, we will see that $f$ is an order-embedding mapping. Given $k_1,k_2\in A$ with $k_1\neq k_2$, the tuples $(0,\dots,0,b\leftarrow_P a_{k_1}^+,0,\dots,0)$ and $(0,\dots,0,b\leftarrow_P a_{k_2}^+,0,\dots,0)$ are clearly different.

It remains to prove that $f$ is an onto mapping. Given a minimal solution $(x_1,\dots,x_m)\in B$, we will prove that there exists an element $k\in A$ such that $f(k)=(x_1,\dots,x_m)$. Taking into account that Equation~\eqref{eq:adjprod2} is solvable and  $a_j^-\neq b$ for each $j\in\{1,\dots,m\}$, by Theorem~\ref{th:adj}, we can ensure that $\max\{a_j^+\mid j\in\{1,\dots,m\}\}\geq b$. In addition, due to $(x_1,\dots,x_m)$ is a solution of Equation~\eqref{eq:adjprod2}, following a similar reasoning to the proof of Theorem~\ref{th:adj}, we obtain that there exists $k\in\{1,\dots,m\}$ such that $a_k^+\geq b$ and $x_k=b\leftarrow_P a_k^+$.

Now, suppose that there exists $j\in\{1,\dots,m\}$ with $j\neq k$ such that $a_j^->b$ and we will obtain a contradiction. Due to $(x_1,\dots,x_m)$ is a solution of Equation~\eqref{eq:adjprod2} and $a_j^->b$, then we can guarantee that $x_j\neq 0$. As product t-norm is a strictly order-preserving mapping, the tuple $(x_1,\dots,x_k,\dots,\frac{x_j}{2},\dots,x_m)$ is also a solution of Equation~\eqref{eq:adjprod2} verifying that $(x_1,\dots,x_k,\dots,\frac{x_j}{2},\dots,x_m)<(x_1,\dots,x_m)$. We have supposed, without loss of generality, that $k<j$. This fact leads us to a contradiction since $(x_1,\dots,x_m)$ is a minimal solution of Equation~\eqref{eq:adjprod2}, by hypothesis. Hence, we can conclude that the inequality $a_j^-\leq b$ holds, for each $j\in\{1,\dots,m\}$ with $j\neq k$. That is, $k\in A$.

Hereafter, we will see that $f(k)=(x_1,\dots,x_m)$. We have already proved that there exists $k\in\{1,\dots,m\}$ such that $x_k=b\leftarrow_P a_k^+$ and the inequality $a_j^-\leq b$ holds, for each $j\in\{1,\dots,m\}$ with $j\neq k$. Suppose that there exists $j\in\{1,\dots,m\}$ with $j\neq k$ such that $x_j>0$ and we will see that this fact contradicts that $(x_1,\dots,x_m)$ is a minimal solution. Since $a_j^-\leq b$, we obtain $a_j^+*0=0$ and $a_j^-*n_P(0)=a_j^-\leq b$. As a result, we can ensure that the tuple $(x_1,\dots,x_{j-1},0,x_{j+1},\dots,x_m)$ is also a solution of Equation~\eqref{eq:adjprod2} satisfying that $(x_1,\dots,x_{j-1},0,x_{j+1},\dots,x_m)<(x_1,\dots,x_j,\dots,x_m)$. This is a contradiction with respect to $(x_1,\dots,x_m)$ is a minimal solution of Equation~\eqref{eq:adjprod2}. Therefore, the equality  $x_j=0$ is verified, for each $j\in\{1,\dots,m\}$ with $j\neq k$. Finally, we can conclude that $f(k)=(0,\dots,0,b\leftarrow_P a_k^+,0,\dots,0)=(x_1,\dots,x_m)$, and thus Statement (3) holds.
\qed
\end{proof}

The following corollary is straightforwardly obtained from Theorem~\ref{th:adj_min}.

\begin{corollary}
Given $a_j^+,a_j^-\in[0,1]$, $b\in(0,1]$, $x_j$ an unknown variable belonging to  $[0,1]$, for each $j\in\{1,\dots,m\}$, and a solvable bipolar max-product FRE as in Equation~\eqref{eq:adjprod2}.
Then, the following statements hold:
\begin{enumerate}
\item[(1)] If there exists $k\in\{1,\dots,m\}$ such that $a_k^-=b$ and $a_j^-\leq b$, for each $j\in\{1,\dots,m\}$, then the least solution of Equation~\eqref{eq:adjprod2} is $(0,\dots,0)$.
\item [(2)] If $a_j^-\neq b$ for each $j\in\{1,\dots,m\}$, then the set of minimal solutions of Equation~\eqref{eq:adjprod2} is given by:
    $$\{(0,\dots,0, x_k,0,\dots,0)\mid  x_k=b\leftarrow_P a_k^+ \hbox{ with } k\in K^+_P\}$$
    where  $K^+_P=\{k\in\{1,\dots,m\}\mid a_k^+\geq b\text{ and }a_j^- < b \text{ for each }j\neq k\}$.
\end{enumerate}
\end{corollary}

Next example clarifies the previous results about the existence both the greatest/least solution and maximal/minimal solutions analyzing the bipolar max-product  FRE given in Example~\ref{ex:adj}.

\begin{example}\label{ex:maximales_minimales}
In Example~\ref{ex:adj} we showed that Equation~\eqref{eq:ex:adj}:

{ \footnotesize \begin{equation*}
(0.4 \ast x_1) \vee (0.7\ast n_P (x_1))\vee(0.2 \ast x_2) \vee (0.1\ast n_P (x_2))\vee(0.5 \ast x_3) \vee (0.2\ast n_P (x_3))=0.3
 \end{equation*} }

is solvable. Now, by Theorems~\ref{th:adj_max} and~\ref{th:adj_min}, we will study whether this equation has a greatest/least solution and maximal/minimal solutions. Taking into account that the hypothesis required in Statement (1) of Theorem~\ref{th:adj_max} is verified, that is $\max\{0.4,0.2,0.5\}\geq0.3$, we can guarantee that the greatest solution exists, which is $(0.75,1,0.6)$.

Applying Theorem~\ref{th:adj_min} (since condition required in Statement (3) is verified) we obtain that this equation has only one minimal solution since:
\[card(\{k\in\{1,2,3\}\mid a_k^+\geq b\text{ and }a_j^- < b \text{ for each }j\neq k\})=card\{1\}=1\]

Making the corresponding computations, we have that $(0.75,0,0)$ is a minimal solution. It is worth highlighting that the tuple $(0.75,0,0)$ is not the least solution of Equation~\eqref{eq:ex:adj}. Indeed, we can easily see that the tuples $(x,0,0.6)$ with $0<x<0.75$ are solutions of Equation~\eqref{eq:ex:adj}, being $(0.75,0,0)$ and $(x,0,0.6)$ incomparable tuples. Moreover, as we show below, the tuple $(0,0,0.6)$ is not a solution of Equation~\eqref{eq:ex:adj}:
 \[(0.4 \ast 0) \vee (0.7\ast 1)\vee(0.2 \ast 0) \vee (0.1\ast 1)\vee(0.5 \ast 0.6) \vee (0.2\ast 0)=0.7\neq0.3\]

Now, modifying the previous equation, we can consider Equation~\eqref{eq:ex:adj_maxmin3} and apply Theorems~\ref{th:adj} and~\ref{th:adj_min} in order to assert that Equation~\eqref{eq:ex:adj_maxmin3} is solvable and it does not have minimal solutions.
{\footnotesize
 \begin{equation}\label{eq:ex:adj_maxmin3}
 (0.4 \ast x_1) \vee (0.7\ast n_P (x_1))\vee(0.2 \ast x_2) \vee (0.1\ast n_P (x_2))\vee(0.5 \ast x_3) \vee (0.9\ast n_P (x_3))=0.3
 \end{equation}
 }
This fact is due to the number of minimal solutions in Equation~\eqref{eq:ex:adj_maxmin3} is:
\[card(\{k\in\{1,2,3\}\mid a_k^+\geq b\text{ and } a_j^- < b \text{ for each }j\neq k\})=card\{\varnothing\}=0\]
\qed
\end{example}

Once the theoretical foundations of the solvability of bipolar max-product FREs have been presented, a toy example is shown below in order to illustrate how a bipolar max-product FRE can be used to represent a real-world situation. Specifically, the overheating of a motor will be modeled by means of a bipolar max-product FRE. Furthermore, given an overheating value, the bipolar FRE is used to determine what is the overheating due to and how a technician should perform in order to solve it.

\begin{example}
The suitable behaviour of a motor directly depends on maintaining its temperature under a certain threshold. To carry out this task, a group of experts stated that controlling the level of water, the level of oil and the radiator fan is crucial. 

Consider the variables $x_1,x_2\in[0,1]$ which represent the truth value of \textit{low water} and \textit{low oil}, respectively, where $1$ indicates empty water/oil container and $0$ that the quantity of water/oil has exceeded the permitted limit. Let $x_3\in\{0,1\}$ be a variable such that $x_3=0$ corresponds to the radiator fan is working and $x_3=1$ to  the radiator fan is stopped. The level of overheating of the motor is represented by a value $b\in[0,1]$, where $0$ indicates a correct temperature and $1$ a critical level of overheating.

After a technical study, the experts reached the following conclusions on the performance of the motor:

\begin{itemize}
\item The motor overheating is directly proportional to the lack of water, with proportionality constant $0.4$. Nevertheless, an excess of water turns out to be even more dangerous for the motor, since in such case it overheats at $0.7$.

Hence, the overheating being caused by the level of water can be modeled by means of the expression
$(0.4 \ast x_1) \vee (0.7\ast n_P (x_1))$. Observe that, if the water container is almost full but not exceeding the limits then the level of overheating is low, since $x_1>0$ and therefore $n_P (x_1)=0$.

\item Similarly, the motor also overheats directly proportional to the lack of oil, but in this case the proportionality constant is $0.2$. Besides, if the oil exceeds the permitted maximum, it causes an overheating of $0.1$. This behaviour can be interpreted by using the expression  $(0.2 \ast x_2) \vee (0.1\ast n_P (x_2))$.

\item The radiator fan has a problem and it sometimes suddenly stops. When this happens, the motor overheats up to $0.5$. Furthermore, the usual behaviour of the radiator makes that the motor overheats at $0.2$. In this case, we obtain the expression 
$(0.5 \ast x_3) \vee (0.2\ast n_P (x_3))$.
\end{itemize}

Now, considering that the level of water, the level of oil and the performance of the radiator fan are the only factors which affect to the ovearheating and basing on the experts knowledge, we can model the behaviour of the level of overheating of the motor through the following bipolar max-product fuzzy relation equation:

{ \footnotesize \begin{equation} \label{eq:ejemplo_b}
(0.4 \ast x_1) \vee (0.7\ast n_P (x_1))\vee(0.2 \ast x_2) \vee (0.1\ast n_P (x_2))\vee(0.5 \ast x_3) \vee (0.2\ast n_P (x_3))=b
 \end{equation} }

This equation is a useful tool in order to know the reasons/causes of a given overheating. For example, if a technician observes that the motor presents an overheating of $0.3$, in order to know the values of the other variables (causes) which imply this level of overheating, we need to solve  Equation~\eqref{eq:ejemplo_b} when  $b=0.3$.

Notice that, the solvability of the obtained equation 
was already 
studied in Examples~\ref{ex:adj} and~\ref{ex:maximales_minimales}. Therefore, we obtain that Equation~\eqref{eq:ejemplo_b} is solvable, where $(0.75,1,0.6)$ is its greatest solution and it has only one minimal solution, the tuple $(0.75,0,0)$.

The most critical cases are then $(0.75,1,0.6)$ and $(0.75,0,0)$. First of all, since the tuple $(0.75,0,0)$ is the unique minimal solution of Equation~\eqref{eq:ejemplo_b}, we deduce that there are no solutions such that $x_1=0$, and thus, the water container is not over the permitted limit. Secondly, we can assert that the level of oil is not giving rise to the overheating, as the variable $x_2$ can be either $0$, $1$, or any value among them. Indeed, a careful sight leads us to conclude that the level of oil can be ignored in this matter, since in the worst case it may result in an overheating of $0.2$, but this overheating is already obtained due to the bad conditions of the radiator fan.

Last but not least, from the greatest solution of Equation~\eqref{eq:ejemplo_b}, the tuple $(0.75,1,0.6)$, we deduce that the variable $x_3$ cannot be equal to $1$. Hence, according to the fact that $x_3\in\{0,1\}$, we conclude that $x_3=0$, and thus the radiator fan is working properly.

Consequently, basing on the conclusions which we have obtained, we would suggest to the technician refilling the water container, but being careful so that it does not exceed the limits.
\qed
\end{example}

\subsection{Solving bipolar max-product FRE systems}

After characterizing the solvability of bipolar max-product FREs with different variables and the product negation, and providing information about the algebraic structure of the set of solutions, we will give a further step to our research. The following goal will be to solve bipolar max-product FRE systems.

According to the results presented until now, one can think that the conditions required to solve an arbitrary system of bipolar max-product FREs will be very similar to the ones given in Theorems~\ref{th:adj_simple} and~\ref{th:adj}. However, the following theorem shows that this reasoning is not true.\\

\begin{theorem}\label{th:sys_mxn}
  Let $a_{ij}^+,a_{ij}^-,x_j\in[0,1]$ and $b_i\in(0,1]$, for each $i\in\{1,\dots,n\}$ and $j\in\{1,\dots,m\}$. Then the bipolar max-product  FRE system given by
  \begin{equation}\label{sys:mxn}
      \bigvee_{j=1}^m(a_{ij}^+ \ast x_j) \vee (a_{ij}^-\ast n_P (x_j))=b_i,\qquad i\in\{1,\dots,n\}
  \end{equation}
  is solvable if and only if there exist two index sets $J^+,J^-\subseteq\{1,\dots,m\}$ with $J^+\cap J^-=\varnothing$ such that, one of the following statements is verified, for each $i\in\{1,\dots,n\}$:
  \begin{enumerate}
    \item[(1)] there exists $j\in J^+$ such that $a_{ij}^+\geq b_i$ and $b_i\leftarrow_P a_{ij}^+\leq b_h\leftarrow_P a_{hj}^+$, for each $h\in\{1,\dots,n\}$.
    \item[(2)] there exists $j\in J^-$ such that $a_{ij}^-=b_i$ and $a_{hj}^-\leq b_h$ for each $h\in\{1,\dots,n\}$.
  \end{enumerate}
\end{theorem}
\begin{proof}
Suppose that there exist two index sets $J^+, J^-\subseteq\{1,\dots,m\}$ with $J^+\cap J^-=\varnothing$, such that each $i\in\{1,\dots,n\}$ satisfies either Statement (1) or (2). Consider the set $J_s^-=\{j\in J^-\mid a_{hj}^-\leq b_h\text{, for each }h\in\{1,\dots,n\}\}$ and the tuple $(\hat{x}_1,\dots,\hat{x}_m)$ defined, for each $j\in\{1,\dots,m\}$, as:
  \[\hat{x}_j= \left\{\begin{array}{lc}
             0 & \hbox{ if}\quad j\in J_s^- \\
             \\ \min\{b_k\leftarrow_P a_{kj}^+\mid k\in\{1,\dots,n\}\} & \hbox{otherwise}
             \end{array}
  \right.\]
Now, we will see that $(\hat{x}_1,\dots,\hat{x}_m)$ is a solution of System~\eqref{sys:mxn}. Given $i\in\{1,\dots,n\}$, by hypothesis, one of the following statements is verified:
  \begin{itemize}
    \item there exists $j\in J^+$ such that $a_{ij}^+\geq b_i$ and $b_i\leftarrow_P a_{ij}^+\leq b_h\leftarrow_P a_{hj}^+$, for each $h\in\{1,\dots,n\}$. In this case, $\hat{x}_j=b_i\leftarrow_P a_{ij}^+$ (note that $j\notin J_s^-$ because $j\in J^+$ and $J^+\cap J^-=\varnothing$). 
        By the definition of the product implication, we obtain that $a_{ij}^+*\hat{x}_j=b_i$. Since by hypothesis $b_i> 0$, we can ensure that $\hat{x}_j>0$, and thus $n_P(\hat{x}_j)=0$ and $a_{ij}^-*n_P(\hat{x}_j)=0$.
    \item there exists $j\in J^-$ such that $a_{ij}^-=b_i$ and $a_{hj}^-\leq b_h$ for each $h\in\{1,\dots,n\}$. Notice that, it is straightforwardly obtained $j\in J_s^-$. Therefore, $\hat{x}_j=0$, and clearly $a_{ij}^+*\hat{x}_j=0$ and $a_{ij}^-*n_P(\hat{x}_j)=a_{ij}^-=b_i$.
  \end{itemize}
  Therefore, we can ensure that there exists $j\in\{1,\dots,m\}$ verifying that $(a_{ij}^+*\hat{x}_j)\vee(a_{ij}^-*n_P(\hat{x}_j))=b_i$.

  In the following, we will prove that the elements $\hat{x}_{j'}$ in $(\hat{x}_1,\dots,\hat{x}_m)$ with $j'\in\{1,\dots,m\}$ and $j'\neq j$ verify $(a_{i{j'}}^+*\hat{x}_{j'})\vee(a_{i{j'}}^-*n_P(\hat{x}_{j'}))\leq b_i$.
  We have that, for every $j'\in\{1,\dots,m\}$, the following statements hold:
  \begin{itemize}
    \item If $j'\in J_s^-$ then $a_{hj'}^-\leq b_h$, for each $h\in\{1,\dots,n\}$. In particular, $a_{ij'}^-\leq b_i$. Therefore, as $\hat{x}_{j'}=0$, we obtain that
        \[(a_{i{j'}}^+*\hat{x}_{j'})\vee(a_{i{j'}}^-*n_P(\hat{x}_{j'}))=0\vee a_{i{j'}}^-=a_{i{j'}}^-\leq b_i\]
    \item If $j'\notin J_s^-$, then $\hat{x}_{j'}=\min\{b_h\leftarrow_P a_{hj'}^+\mid h\in\{1,\dots,n\}\}$. On the one hand, as $b_h>0$ for each $h\in\{1,\dots,n\}$, then $b_h\leftarrow_P a_{hj'}^+>0$ for each $h\in\{1,\dots,n\}$. As a consequence, we obtain that $\hat{x}_{j'}>0$, and thus $n_P(\hat{x}_{j'})=0$. On the other hand, by definition of $\hat{x}_{j'}$, we can ensure that, in particular, $\hat{x}_{j'}\leq b_i\leftarrow_P a_{ij'}^+$. Therefore, applying the adjoint property, the next inequality $a_{i{j'}}^+*\hat{x}_{j'}\leq b_i$ holds. As a result,   the following inequality is verified:
 \[(a_{i{j'}}^+*\hat{x}_{j'})\vee(a_{i{j'}}^-*n_P(\hat{x}_{j'}))= (a_{i{j'}}^+*\hat{x}_{j'}) \leq b_i\]
    
    \end{itemize}
    Therefore,
    \[\bigvee_{j=1}^m(a_{ij}^+ \ast \hat{x}_{j}) \vee (a_{ij}^-\ast n_P (\hat{x}_{j}))=b_i\]
By an analogous reasoning for each $i\in\{1,\dots,n\}$, we can conclude that $(\hat{x}_{1},\dots,\hat{x}_{m})$ is a solution of System~\eqref{sys:mxn}.

In order to prove the counterpart, suppose that System~\eqref{sys:mxn} is solvable, and let build two sets $J^+,J^-\subseteq\{1,\dots,m\}$ with $J^+\cap J^-=\varnothing$ such that, for each $i\in\{1,\dots,n\}$, either Statement $(1)$ or Statement $(2)$ is verified.

Given a solution $(\hat{x}_{1},\dots,\hat{x}_{m})$ of System~\eqref{sys:mxn} and two index sets $J^+, J^-$ defined as follows:
\begin{eqnarray*}
J^+&=&\{j\in\{1,\dots,m\}\mid \hat{x}_j>0\}\\
J^-&=&\{j\in\{1,\dots,m\}\mid \hat{x}_j=0\}
\end{eqnarray*}
Clearly, $J^+\cap J^-=\varnothing$. Given $i\in\{1,\dots,n\}$, as $(\hat{x}_{1},\dots,\hat{x}_{m})$ is a solution of System~\eqref{sys:mxn}, one of the following statements is verified:
\begin{itemize}
      \item[(a)] there exists $j\in\{1,\dots,m\}$ such that $a_{ij}^+*\hat{x}_{j}=b_i$.
      \item[(b)] there exists $j\in\{1,\dots,m\}$ such that $a_{ij}^-*n_P(\hat{x}_{j})=b_i$.
\end{itemize}

    On the one hand, if Statement (a) holds, as $b_i>0$, we deduce that $\hat{x}_{j}>0$, and thus $j\in J^+$. Moreover, since operator $*$ is monotonic and $a_{ij}^+*1=a_{ij}^+$, the equality $a_{ij}^+*\hat{x}_{j}=b_i$ implies that $a_{ij}^+\geq b_i$. In fact, by a similar reasoning to the one given in Theorem~\ref{th:adj_simple}, it can be proved that  $\hat{x}_{j}=b_i\leftarrow_P a_{ij}^+$.

    Finally, suppose that there exists $h\in\{1,\dots,n\}$ such that $b_i\leftarrow_P a_{ij}^+> b_h\leftarrow_P a_{hj}^+$, this is, $\hat{x}_{j}> b_h\leftarrow_P a_{hj}^+$. From the adjoint property, we obtain that $\hat{x}_{j}\leq b_h\leftarrow_P a_{hj}^+$ if and only if $a_{hj}^+*\hat{x}_{j}\leq b_h$, Therefore, we can ensure that $a_{hj}^+*\hat{x}_{j}>b_h$. This fact leads us to a contradiction, since $(\hat{x}_{1},\dots,\hat{x}_{m})$ would not be a solution of equation $h$ in System~\eqref{sys:mxn} and so, it would not be  a solution of the mentioned system.

    Hence, we can ensure that, if Statement (a) holds, then there exists $j\in J^+$ such that $a_{ij}^+\geq b_i$ and $b_i\leftarrow_P a_{ij}^+\leq b_h\leftarrow_P a_{hj}^+$, for each $h\in\{1,\dots,n\}$. That is, Statement (1) holds.

    On the other hand, suppose that Statement (b) is verified. Since $b_i> 0$, the equality $a_{ij}^-*n_P(\hat{x}_{j})=b_i$ implies that $n_P(\hat{x}_{j})>0$, and thus $\hat{x}_{j}=0$. That is, $j\in J^-$. Consequently, we obtain that $n_P(\hat{x}_{j})=1$, and straightforwardly $a_{ij}^-=b_i$.

   To end, suppose that there exists $h\in\{1,\dots,n\}$ such that $a_{hj}^->b_h$. 
   As a result, we obtain that $(\hat{x}_{1},\dots,\hat{x}_{m})$ is not a solution of equation $h$ in System~\eqref{sys:mxn} and thus it is not a solution of that system. That is, we obtain a contradiction.

    Therefore, if Statement (b) holds, then there exists $j\in J^-$ such that $a_{ij}^-=b_i$ and $a_{hj}^-\leq b_h$ for each $h\in\{1,\dots,n\}$. In other words, Statement (2) holds.
    \qed
\end{proof}

The following example clarifies the previous result.

\begin{example}
Given the following bipolar max-product FRE system with three equations and three unknown variables

{ \footnotesize
\begin{eqnarray*}
(1*x_1) \vee (0.7*n_P(x_1)) \vee (0.2*x_2) \vee (0.4*n_P(x_2)) \vee (0.5*x_3) \vee (0.7*n_P(x_3))\!\!\!\!&=&\!\!\!\!0.7\\
(0.8*x_1) \vee (0.1*n_P(x_1)) \vee (0.8*x_2) \vee (0.2*n_P(x_2)) \vee (0.3*x_3) \vee (0.6*n_P(x_3))\!\!\!\!&=&\!\!\!\! 0.6\\
(0.4*x_1) \vee (0.2*n_P(x_1)) \vee (0.3*x_2) \vee (0.3*n_P(x_2)) \vee (0.6*x_3) \vee (0.2*n_P(x_3))\!\!\!\!&=&\!\!\!\! 0.3
\end{eqnarray*}
}

We will consider the index sets $J^+=\{2,3\}$ and $J^-=\{1\}$. We will see that for each $i\in\{1,2,3\}$ the conditions required in the hypothesis of Theorem~\ref{th:sys_mxn} are verified. First of all, notice that $J^+\cap J^-=\varnothing$.
\begin{itemize}
  \item Case $i=1$ (first equation of the system):  there exists $j=1$ belonging to $J^-$ satisfying that $a_{11}^-=b_1=0.7$, $a_{21}^-=0.1\leq0.6= b_2$ and $a_{31}^-=0.2\leq 0.3=b_3$.
  \item Case $i=2$ (second equation of the system): we obtain that the index $j=2\in J^+$ verifies that $0.8=a_{22}^+\geq b_2=0.6$. In addition, the inequalities $b_2\leftarrow_P a_{22}^+\leq b_3\leftarrow_P a_{32}^+$ and $b_2\leftarrow_P a_{22}^+\leq b_1\leftarrow_P a_{12}^+$ hold, as we see below:
      \[0.6\leftarrow_P 0.8=0.75\leq 1=0.7\leftarrow_P 0.2\]
      \[0.6\leftarrow_P 0.8=0.75\leq 1=0.3\leftarrow_P 0.3\]
  \item Case $i=3$ (third equation of the system):  we find the index $j=3$ belonging to $J^+$ satisfying that $0.6=a_{33}^+\geq b_3=0.3$. Furthermore, the next inequalities are obtained:
      \[b_3\leftarrow_P a_{33}^+=0.3\leftarrow_P 0.6=0.5\leq 1=0.7\leftarrow_P 0.5=b_1\leftarrow_P a_{13}^+\]
      \[b_3\leftarrow_P a_{33}^+=0.3\leftarrow_P 0.6=0.5\leq 1=0.6\leftarrow_P 0.3=b_2\leftarrow_P a_{23}^+\]
\end{itemize}

Consequently, Theorem~\ref{th:sys_mxn} can be applied which leads us to obtain that the given system is solvable. Indeed, from the computations above, one can easily check that the tuple $(0,0.75,0.5)$ is a solution of such system.

Observe that, the tuple $(0,0.75,0.5)$ is not the unique solution of the system. We can find some other different definitions of the index sets $J^+$ and $J^-$ such that, for each $i\in\{1,2,3\}$ the hypothesis of Theorem~\ref{th:sys_mxn} are satisfied. For instance, considering the index sets $J^+=\varnothing$ and $J^-=\{1,2,3\}$, it is easy to deduce that the system is solvable and the tuple $(0,0,0)$ is a solution of the given system.

It is also important to mention that,  the set $J^+$ does not need to be the complementary set of $J^-$. That is, the equality $J^+\cup J^-=\{1,2,3\}$ is not required in order to guarantee the solvability of the system. For example, the index sets $J^+=\{1\}$ and $J^-=\{3\}$ verify the hypothesis of Theorem~\ref{th:sys_mxn} and clearly $J^+\cup J^-=\{1,3\}\neq \{1,2,3\}$. From the sets $J^+=\{1\}$ and $J^-=\{3\}$, we obtain the solution $(0.75,0,0)$.
\qed
\end{example}

We have provided sufficient and necessary conditions which allow us to know when an arbitrary bipolar max-product FRE system with the residuated negation of the product t-norm is solvable. In the future, we will study the resolution of bipolar max-product FRE systems in which independent terms can take the value zero.

\section{Conclusions and future work}
A broad study on the resolution of bipolar max-product FREs has been carried out, considering the non-involutive negation operators defined from the implication associated with the product t-norm. Three different parts can be distinguished in our research according to the complexity of the considered bipolar equations. The first part introduces a characterization on the resolution of bipolar max-product FREs with only one unknown variable. The second part shows under what conditions a bipolar max-product fuzzy relation equation containing different variables is solvable. Moreover, interesting properties related to the algebraic structure of the set of solutions have been included. The third part considers bipolar max-product FRE systems and presents the requirements to guarantee when these systems are solvable. Notice that, the residuated negations related to the G\"odel implication and the product implication coincide, and therefore, we can ensure that the solvability for bipolar max-product FREs with the G\"odel negation has also been analyzed in this paper.

As future work, based on Theorem~\ref{th:sys_mxn}, we will study the introduction of an efficient algorithm,   determining  whether a system of bipolar max-product fuzzy relation equations is solvable and   computing at least one solution. It will be also fundamental to investigate more properties in order to know the algebraic structure of the complete set of solutions corresponding to an arbitrary solvable system of bipolar max-product fuzzy relation equations. 
Moreover, we are interested in studying bipolar fuzzy relation equations based on other general operators, such as uninorms, u-norms, adjoint triples, etc.

{

}

\end{document}